# MULTISCALE LIKELIHOOD ANALYSIS AND COMPLEXITY PENALIZED ESTIMATION[1]

By Eric D. Kolaczyk and Robert D. Nowak

*Boston University and University of Wisconsin*

We describe here a framework for a certain class of multiscale likelihood factorizations wherein, in analogy to a wavelet decomposition of an $L^2$ function, a given likelihood function has an alternative representation as a product of conditional densities reflecting information in both the data and the parameter vector localized in position and scale. The framework is developed as a set of sufficient conditions for the existence of such factorizations, formulated in analogy to those underlying a standard multiresolution analysis for wavelets, and hence can be viewed as a multiresolution analysis for likelihoods. We then consider the use of these factorizations in the task of nonparametric, complexity penalized likelihood estimation. We study the risk properties of certain thresholding and partitioning estimators, and demonstrate their adaptivity and near-optimality, in a minimax sense over a broad range of function spaces, based on squared Hellinger distance as a loss function. In particular, our results provide an illustration of how properties of classical wavelet-based estimators can be obtained in a single, unified framework that includes models for continuous, count and categorical data types.

**1. Introduction.** Wavelet-based methods have had a decided impact on the field of nonparametric function estimation in the past decade, particularly where concerned with inhomogeneous objects, as might be encountered in applications such as signal and image processing. The near-optimality of their risk properties (in a minimax sense) and their adaptivity to various ranges of unknown degrees of smoothness, combined with simple and efficient algorithms for practical implementation, have all contributed to this

Received August 2001; revised August 2002.

[1]Supported by NSF Grants BCS-00-79077 and MIP-97-01692, Office of Naval Research Grants N00014-99-1-0219 and N00014-00-1-0390 and Army Research Office Grant DAAD19-99-1-0290.

impact. See Donoho, Johnstone, Kerkyacharian and Picard (1995), for example, and the discussions therein.

Much of this work rests upon a framework that assumes a standard Gaussian "signal-plus-noise model," that is, $X_i = f(i/N) + Z_i$, where the $Z_i$ are i.i.d. standard normal random variables and the $f(i/N)$ are equispaced samples of an unknown function $f$ on the unit interval. This model is then combined with an expansion

$$f(t) = \sum_{(j,k) \in \mathbb{Z}^2} \omega_{j,k} \psi_{j,k}(t), \tag{1}$$

through which the details in $f$ gained between certain approximations at scales indexed by $j$ and $j+1$ (formally, scales $2^{-j}$ and $2^{-(j+1)}$), in the vicinity of locations indexed by $k$, are captured by the coefficients $\omega_{j,k} \equiv \langle f, \psi_{j,k} \rangle$. The $\psi_{j,k}(t) \equiv 2^{j/2} \psi(2^j t - k)$ are orthonormal dilations–translations of a wavelet function $\psi(t)$ satisfying the admissibility condition $\int \psi(t)\, dt = 0$, as well as various other conditions on smoothness, symmetry or such as desired.

Little or no work, however, has been done extending the wavelet paradigm to certain other common noise models. We have in mind, in particular, models for count and categorical data, such as Poisson or multinomial. Count data of this sort arises in a variety of contexts, such as high-energy astrophysics or medical imaging, while a good example of such categorical data might be the images found in landcover classification from remote sensing data. The authors of this paper, in addition to various collaborators, have in recent years pursued a program that seeks to extend wavelet-based frameworks in such directions through the use of various multiscale probability models [e.g., Kolaczyk (1999a), Timmermann and Nowak (1999), Nowak and Kolaczyk (2000) and Kolaczyk and Huang (2001)]. At the heart of this program is the concept of a multiscale factorization of a given data likelihood, $p(\mathbf{X}|\boldsymbol{\theta})$, in analogy to the orthogonal wavelet decomposition in (1), that is, expressions like

$$p(\mathbf{X}|\boldsymbol{\theta}) \propto \prod_{j,k} p(X_{j+1,2k}|X_{j,k}, \omega_{j,k}), \tag{2}$$

where $X_{j,k}$ contains information in the original data $\mathbf{X}$ local to scale $j$ and position $k$, $X_{j+1,2k}$ contains information within a refined subregion of that and the parameters $\omega_{j,k}$ reflect similar information in the original parameter $\boldsymbol{\theta}$. Note that the pursuit of such factorizations differs from attempting to simply work with the likelihood induced by applying a wavelet transform to the original data, as the latter approach tends quickly to lead to model expressions that suffer from difficulty of interpretation and computational intractability outside of the Gaussian case [Kolaczyk (1999b)].

Our goal in this paper is to show that a systematic approach can be taken to the topic of multiscale probability models, in which for a certain class of such models the relevant characteristics of traditional wavelet-based models and their extensions are paralleled quite closely. In particular, our contribution consists of two related components. First, we show that factorizations like that in (2) arise when conditions for a "multiresolution analysis (MRA)" of likelihoods are satisfied, where the conditions are a blending of concepts from the fields of wavelets, recursive partitioning and graphical models. These conditions are then shown to characterize the Gaussian, Poisson and multinomial models. Hence, such multiscale probability models provide an example of a unified framework for modeling data of continuous, count and categorical types in a fashion sensitive to location–scale variation. Second, we quantify the risk behavior of certain nonparametric, complexity penalized likelihood estimators based on our factorizations. We show that a near-optimality and adaptivity completely analogous to that of wavelet-based estimators holds for these disparate data types, for appropriately defined smoothness classes, using the squared Hellinger distance as a loss function. The technical details behind these results rely on upper bounds on the risk in the spirit of recent work by Birgé and Massart [e.g., Barron, Birgé and Massart (1999) and references therein], but these bounds are derived by adapting a technique of Li (1999) and Li and Barron (2000). In addition to the above two primary contributions, we also comment briefly on the algorithmic efficiency with which our various estimators may be calculated, an indication of their relevance to practice as well as theory.

The body of the paper is arranged as follows. In Section 2 we present our multiresolution analysis for likelihoods. Then in Section 3 we provide necessary details for a certain class of models for continuous, count and categorical data types, and introduce three estimators of the relevant underlying function. Risk properties of these estimators are then stated in Section 4. Proofs of these results are detailed in Section 5, and some final comments and discussion are compiled in Section 6. Finally, a result on the algorithmic complexity of our estimators is proven in the Appendix.

**2. A multiresolution analysis for likelihoods.** Consider a stochastic process $X(t)$ on the interval $[0,1)$ that, either by choice or perhaps by the limitations of measuring instruments, is observed only discretely on the intervals $I_i \equiv [i/N, (i+1)/N)$, $i = 0, \ldots, N-1$. Furthermore, suppose that corresponding to this process is a function $\theta(t), t \in [0,1)$. We will assume that the effect of the discretization is to yield a vector of measurements $\mathbf{X} \equiv (X_0, \ldots, X_{N-1})$, associated with a vector of parameters $\boldsymbol{\theta} \equiv (\theta_0, \ldots, \theta_{N-1})$, where each pair $(X_i, \theta_i)$ corresponds to the interval $I_i$ and is obtained by sampling the function $\theta(\cdot)$ and then sampling $X_i$, in a manner to be made precise later. We will denote the likelihood of $\mathbf{X}$, given the parameter value $\boldsymbol{\theta}$, by $p(\mathbf{X}|\boldsymbol{\theta})$ generically for both discrete and continuous distributions on $\mathbf{X}$ [i.e., $p(\cdot)$ is to be understood to be defined with respect to an appropriate measure $\nu$]. At times, when convenient, we may abbreviate this notation as $p_{\boldsymbol{\theta}}$.

Informally speaking, a simple yet standard multiscale analysis of the data $\mathbf{X}$ is achieved by defining the dyadic intervals $I_{j,k} \equiv [k/2^j, (k+1)/2^j)$, for $j = 0, \ldots, J-1$, $k = 0, \ldots, 2^j - 1$, and $J = \log_2(N)$ (i.e., with $N$ assumed to be a power of 2, for convenience), and associating with each a summary statistic $X_{j,k} \equiv \sum_{i/N \in I_{j,k}} X_i$. That is, we define an analysis separating the information in $\mathbf{X}$ into its components at various combinations of scale and position $(j, k)$. This strategy, of course, underlies the analysis of $\mathbf{X}$ with respect to an orthonormal basis of dyadic Haar wavelets, specifically, analysis through the discrete inner products of $\mathbf{X}$ with functions $h_{j,k}(i) \equiv (\chi_{j+1,2k+1}(i) - \chi_{j+1,2k}(i))/N_{j,k}^{1/2}$, defined on the index set $\{0, 1, \ldots, N-1\}$, where $\chi_{j,k}$ is the characteristic function for the discrete analogue of the interval $I_{j,k}$, that is, $\{i : I_i \subseteq I_{j,k}\}$, and $N_{j,k}$ is the cardinality of this set.

This particular notion of multiscale analysis can be generalized by generalizing the underlying notion of partitioning. Specifically, beginning with the unit interval $[0, 1)$, we partition that interval in a recursive fashion, where split points are constrained to the endpoints of the original sampling intervals $I_i$, until a complete recursive partition (C-RP) $\mathcal{P}^* \equiv \{I_i\}_{i=0}^{N-1}$ is achieved. That is, beginning with the trivial partition $[0, 1)$, we split that into two pieces at one of the points $\{i/N\}_{i=1}^{N-1}$. Then, proceeding in a recursive fashion, given a partition $\mathcal{P}$ intermediate to $[0, 1)$ and $\mathcal{P}^*$, we refine that partition by splitting one and only one of the intervals $I \in \mathcal{P}$ at one of the remaining allowable points (i.e., those points in the intersection of $\{i/N\}_{i=1}^{N-1}$ and the interior of $I$). We often will call the interval $I$ in such cases the "parent" interval, and the two corresponding subintervals, say $I_{\text{ch}(I),\text{l}}$ and $I_{\text{ch}(I),\text{r}}$, the left and right interval "children." Partitions $\mathcal{P}'$ produced further along in the recursive process than a partition $\mathcal{P}$ will be said to be refinements of $\mathcal{P}$, which we will denote $\mathcal{P} \prec \mathcal{P}'$ (refinement that includes potential equivalence will be denoted using "$\preceq$"). Finally, for a given $\mathcal{P} \preceq \mathcal{P}^*$, let $\mathcal{I}(\mathcal{P})$ be the collection of all intervals $I$ found in at least one partition $\mathcal{P}' \preceq \mathcal{P}$, and let $\mathcal{I}_{\text{NT}}(\mathcal{P})$ be all such nonterminal intervals [i.e., all intervals $I \in \mathcal{I}(\mathcal{P})$ that are *not* in $\mathcal{P}$ itself].

The multiscale analysis of $\mathbf{X}$ corresponding to $\mathcal{P}^*$ is then composed of the statistics $X_I \equiv \sum_{i/N \in I} X_i$, for all intervals $I \in \mathcal{I}(\mathcal{P}^*)$. This analysis can be linked in turn to analysis with respect to an orthonormal basis of so-called unbalanced Haar wavelets [Girardi and Sweldens (1997)]

$$(3) \qquad h_I(i) = c_I' \left[ \frac{\chi_{\text{ch}(I),\text{r}}(i)}{N_{\text{ch}(I),\text{r}}} - \frac{\chi_{\text{ch}(I),\text{l}}(i)}{N_{\text{ch}(I),\text{l}}} \right],$$

where $N_I = \#\{i : I_i \subseteq I\}$ is the discrete length of an interval $I$, and $c'_I = (N^{-1}_{\text{ch}(I),\text{r}} + N^{-1}_{\text{ch}(I),\text{l}})^{-1/2}$ is a normalizing constant. Note that the dyadic analysis above is seen to be the special case in which parent intervals are split only into two interval "children" of equal size, that is, $N_{\text{ch}(I),\text{l}} = N_{\text{ch}(I),\text{r}} = N_I/2$, yielding the complete recursive dyadic partition (C-RDP) $\mathcal{P}^*_{\text{Dy}}$ and the dyadic Haar wavelets $h_{j,k}$.

Our goal in this section is to show how the above concepts may be used to produce a probabilistic analogue of an orthonormal wavelet expansion like that in (1). We do so by introducing a formal analogue of the key conceptual framework underlying the latter, that is, multiresolution analysis.

### 2.1. Development of a formal multiresolution analysis.

2.1.1. *Function space multiresolution analysis.* Fundamental to the concept of wavelets is the notion of a multiresolution analysis (MRA). Briefly, the idea behind this method is to construct a sequence of subspaces $V_j \subseteq L^2(\mathbb{R})$, across scales $j$, whose members contain successively finer approximations to functions $f \in L^2(\mathbb{R})$. The classical multiresolution analysis [e.g., see Daubechies (1992)] requires the following three sets of characteristics of these subspaces.

(A) Hierarchy of nested subspaces. The subspaces $V_j$ satisfy the condition

$$\cdots V_{-2} \subset V_{-1} \subset V_0 \subset V_1 \subset V_2 \cdots,$$

where $\bigcap_{j \in \mathbb{Z}} V_j = \{0\}$ and $\overline{\bigcup_{j \in \mathbb{Z}} V_j} = L^2(\mathbb{R})$.

(B) Orthonormal basis within $V_0$. There exists a function $\phi$ such that the collection $\{\phi(\cdot - k)\}_{k \in \mathbb{Z}}$ forms an orthonormal basis for $V_0$.

(C) Scaling between and translation within subspaces.

$$\begin{aligned} g \in V_j &\iff g(2^{-j} \cdot) \in V_0, \\ g \in V_0 &\implies g(\cdot - k) \in V_0 \qquad \forall\, k \in \mathbb{Z}. \end{aligned}$$

Our interest in the above characteristics (gathered into these three labeled categories for our own later convenience of exposition) centers primarily on the fact that they form a set of sufficient conditions for the existence of a (wavelet) function $\psi \in L^2(\mathbb{R})$ for which the collection $\{\psi_{j,k}\}$ forms an orthonormal basis of $L^2(\mathbb{R})$, as in (1). In other words, these conditions assure a multiscale decomposition or "decoupling" of any given function $f \in L^2(\mathbb{R})$ into components of $L^2$ "energy" localized to certain combinations of scale $j$ and position $k$. The fact that this decoupling is with respect to an orthonormal basis implies that knowledge of these components (i.e., the coefficients and their corresponding wavelets) is equivalent to knowledge of the function $f$ itself—only the representation has changed.

2.1.2. *Likelihood multiresolution analysis.* In analogy to the three conditions (A)–(C) outlined in Section 2.1.1, we provide four conditions $(A^*)$– $(D^*)$ sufficient to insure a certain multiscale likelihood factorization. The first three conditions will be seen to play roles that parallel those of (A)– (C). However, to obtain a factorization fully analogous to an orthonormal basis decomposition, the fourth condition $(D^*)$ is needed. Our conditions are as follows.

$(A^*)$ Hierarchy of recursive partitions. A hierarchy of recursively defined partitions

$$\cdots \mathcal{P}_{\ell-1} \prec \mathcal{P}_\ell \prec \mathcal{P}_{\ell+1} \cdots,$$

beginning with $[0,1)$ and ending with a C-RP $\mathcal{P}^* = \{I_i\}_{i=0}^{N-1}$.

$(B^*)$ Independence within $\mathcal{P}^*$. The components of $\mathbf{X}$ are statistically independent, the components of $\boldsymbol{\theta}$ are $L$-independent with respect to the likelihood of $\mathbf{X}$, that is,

$$p(\mathbf{X}|\boldsymbol{\theta}) = \prod_{i=0}^{N-1} p(X_i|\theta_i),$$

and the p.d.f. for each $X_i$ is a member of some common parametric family $\mathcal{F} \equiv \{p(\cdot|\theta) : \theta \in \Theta \subseteq \mathbb{R}\}$.

$(C^*)$ Reproducibility between partitions. The family $\mathcal{F}$ is reproducible in $\theta$, in the sense that, for all $I \in \mathcal{I}(\mathcal{P}^*)$ and $\forall \boldsymbol{\theta} \in \Theta^{N_I}$, the p.d.f. of $X_I \equiv \sum_{i/N \in I} X_i$ is $p(X_I|\theta_I) \in \mathcal{F}$, where $\theta_I \equiv \sum_{i/N \in I} \theta_i$.

$(D^*)$ "Decoupling" of parameters with partitions (i.e., cuts). For any $X_i \sim p(\cdot|\theta_i) \in \mathcal{F}, i \in \{i_1, i_2\}$, there exists some reparameterization $(\theta_{i_1}, \theta_{i_2}) \to (\theta, \omega)$ such that

$$p(X_{i_1}, X_{i_2}|\theta_{i_1}, \theta_{i_2}) = p(X|\theta)p(X_{i_1}|X, \omega),$$

where $X \equiv X_{i_1} + X_{i_2}$ and $\theta \equiv \theta_{i_1} + \theta_{i_2}$. That is, the sum $X$ is a cut [e.g., Barndorff-Nielsen (1978)] for $(X_{i_1}, X_{i_2})$.

Some remarks on these conditions are useful prior to stating our main result on likelihood factorizations. First, note that through $(A^*)$ the notion of multiple resolutions takes the form of recursive partitioning (or, conversely, hierarchical aggregation) of our data space. Second, in condition (B$^*$), the assumption of a likelihood factorization in the original data space, with respect to the index set $\{0, \ldots, N-1\}$, and in components identical up to the parameters $\theta_i$, mirrors the spirit and function of the orthonormal basis $\{\phi(\cdot - k)\}$ in $V_0$. The condition of $L$-independence requires that the domain of variation of $\boldsymbol{\theta}$ be equal to the product of the domains of the components $\theta_i$ and that the role of $\boldsymbol{\theta}$ in the likelihood $p(\mathbf{X}|\boldsymbol{\theta})$ can be separated in one-to-one correspondence with the statistically independent likelihood components of the $X_i$ [Barndorff-Nielsen (1978)]. Third, for each interval $I \in \mathcal{I}(\mathcal{P}^*)$, it is desirable to combine the information in $\{X_{i_1}, \ldots, X_{i_{N_I}}\}$ into a single summary statistic, in analogy to orthogonal projection of a function onto subspaces $V_j$. We do so using the simplest approach, direct summation, that is, $X_I = \sum_{i/N \in I} X_i$. Condition (C$^*$) dictates that the distributional family $\mathcal{F}$ is in some sense "invariant" under this summation—a scale-invariance, in a sense. Practically speaking, there are in fact a number of similar definitions available. We use the very simplest definition here, found, for example, in Wilks (1962), which describes the well-known behavior of such distributions as the Gaussian, Poisson, Cauchy and others.

Our perspective in introducing these conditions is one in which we view an orthonormal basis decomposition essentially as a "decoupling" of information over some meaningful index space. In the case of wavelets, the indexing $(j, k)$ refers to information local in scale and position. For likelihoods, a fully analogous decoupling requires both independence and $L$-independence in this indexing, such as we assume holds true in the original indexing $i$ through condition (B$^*$). Conditions (A)–(C) are sufficient to guarantee a multiscale decoupling of a function $f$ in the manner of (1). However, conditions (A$^*$)–(C$^*$) yield only statistical independence in the multiscale indexing, and not $L$-independence. Put another way, we have a factorization of $p(\mathbf{X}|\boldsymbol{\theta})$ into components that are functions of only *local* information in $\mathbf{X}$, but possibly *global* information in $\boldsymbol{\theta}$. Condition (D$^*$) remedies this situation.

We now state the main result of this section.

THEOREM 1. *Assume that the conditions* (A$^*$)–(D$^*$) *hold. Then there exists a factorization of the form*

(4) $$p(\mathbf{X}|\boldsymbol{\theta}) = p(X_{I_{00}}|\theta_{I_{00}}) \prod_{I \in \mathcal{I}_{\mathrm{NT}}(\mathcal{P}^*)} p(X_{\mathrm{ch}(I),\mathrm{l}}|X_I, \omega_I),$$

*with respect to some reparameterization* $\{\theta_{I_{00}}, \boldsymbol{\omega}\}$ *of* $\boldsymbol{\theta}$, *for* $I_{00} \equiv [0, 1)$ *and* $\theta_{I_{00}} \equiv \sum_{i=0}^{N-1} \theta_i$.

Proof of the theorem follows immediately, in light of the conditions and the above discussion. Alternatively, this result may be viewed as a consequence of the fact that conditions (A$^*$) and (B$^*$) imply a so-called directed, local Markov property for the graphical model given by $\{X_I\}_{I \in \mathcal{I}(\mathcal{P}^*)}$, where the underlying graph is just a binary tree $\mathcal{T}^* \equiv \mathcal{T}(\mathcal{P}^*)$ equipped with arrows denoting parent–child relationships. Equation (4) then follows as an example of a recursive factorization [e.g., Lauritzen (1996), Theorem 3.27]. The fact that the conditional distributions in the factorization are of the same family follows from condition (C*), and the reparameterization follows from condition (D*).

For a random variable **X** associated with family $\mathcal{F}$ for which (A*)–(D*) are satisfied, we will say that $\mathcal{F}$ *allows a likelihood MRA with respect to* $\boldsymbol{\theta}$. In considering the factorization in (4), note that the role of a wavelet coefficient–function pair $(\omega_{j,k}, \psi_{j,k})$, in capturing detail lost between scales $j+1$ and $j$ in approximating $f \in L^2(\mathbb{R})$, is played here by the conditional density $p(X_{\mathrm{ch}(I),1}|X_I, \omega_I)$, a natural form of expressing the information lost between the aggregations dictated by a partition $\mathcal{P} \prec \mathcal{P}^*$ and its immediate predecessor.

2.2. *Characterization.* The conditions of Theorem 1 may be used to characterize certain families $\mathcal{F}$ that allow a multiresolution analysis. We illustrate with the canonical case in which $\mathcal{F}$ is a one-parameter natural exponential family (NEF), that is,

$$p(X_i|\eta_i) = a(\eta_i)b(X_i)\exp\{\eta_i X_i\}$$

with respect to some sigma-finite measure $\nu(\cdot)$, for natural parameter $\eta \in \mathrm{E} \subset \mathbb{R}$. Specifically, we have the following result.

THEOREM 2. *Suppose that $\mathcal{F}$ is a (minimal and steep) one-parameter NEF. Then it follows that:*

(i) *$\mathcal{F}$ allows a likelihood MRA with respect to the natural parameterization $\boldsymbol{\theta} \equiv \boldsymbol{\eta}$ if and only if $\mathcal{F}$ is the family of Gaussian distributions;*

(ii) *$\mathcal{F}$ allows a likelihood MRA with respect to the mean parameterization $\boldsymbol{\theta} \equiv \boldsymbol{\mu}(\boldsymbol{\eta})$ if and only if $\mathcal{F}$ is either the family of Gaussian distributions or the family of Poisson distributions.*

Proof of Theorem 2 follows from the use of results in the literature on reproducibility and cuts. One begins by noting that the collection of NEFs $\mathcal{F}$ satisfying (C*) must be contained within the collection of such $\mathcal{F}$ which do so in the case of i.i.d. random variables, that is, where $\eta_0 = \cdots = \eta_{N-1} \equiv \eta$, for some $\eta \in \mathrm{E}$. A characterization of this latter case, under a slight generalization of our own definition of reproducibility, is provided in Bar-Lev and Enis (1986). Specifically, among various other results, these authors show that reproducibility implies that $\mathcal{F}$ must have a power variance function (PVF) and that there are only four NEF-PVF families. Examination of the cumulant generating function for these four families yields

MULTISCALE LIKELIHOOD ANALYSIS 9

the candidate distributions under cases (i) and (ii) of the theorem. The result follows by confirming that the sum of independent random variables forms a cut for the joint distribution in the case of the Gaussian and Poisson families, hence satisfying condition ($D^*$), which is straightforward [e.g., Barndorff-Nielsen (1978)].

Theorem 2 thus establishes formally a role for the Gaussian and Poisson distributions in our class of multiscale probability models. These models have been derived from first principles in previous work [e.g., Kolaczyk (1999a), Timmermann and Nowak (1999) and Nowak (1999)]. Similarly, a moment's thought reveals that the factorization in (4) holds as well for the case in which **X** follows a multinomial distribution, given the appearance of that distribution when conditioning a vector of independent Poisson random variables on their total (i.e., $X_{I_{00}}$). In fact, it can be shown using results from the literature on cuts for discrete NEFs [e.g., Barndorff-Nielsen (1978) or, alternatively, the work of Joshi and Patil (1970)], that (4) holds for all members of the class of sum-symmetric power series distributions (SSPSD), that is, where **X** has probability mass function of the form

$$(5) \qquad p(\mathbf{X}|\boldsymbol{\theta}) = b(\mathbf{X}) \frac{\theta_0^{X_0} \cdots \theta_{N-1}^{X_{N-1}}}{g(\boldsymbol{\theta})},$$

where the generating function $g(\cdot)$ depends on $\boldsymbol{\theta}$ only through $\theta_0 + \cdots + \theta_{N-1}$. The Poisson and multinomial families are two members of this class, with the former being the unique member for which the components of **X** are uncorrelated. Hence, this result also indicates that the statistical independence assumed in condition ($B^*$), while sufficient, is not necessary for the result of Theorem 1.

To what degree these results may be extended remains an open question. The above discussion suggests that extensions are unlikely without significant relaxation of conditions ($A^*$)–($D^*$). While such extensions potentially could be interesting, for example, in the event that they may necessarily parallel certain aspects of the "second generation" [Sweldens (1998)] extensions of the classical MRA underlying conditions (A)–(C), it is not clear whether they would lead to methods of practical interest.

**3. Multiscale penalized maximum likelihood estimation.** We now turn our attention to the problem of estimating the unknown parameter vector $\boldsymbol{\theta}$ from data **X**, when the underlying distributional family allows a likelihood MRA. In light of the results of the previous section, for the remainder of this paper we will restrict our attention to three models, those of the Gaussian, Poisson and multinomial families, as canonical examples of models for continuous, count and categorical measurements. Additionally, in preparation for the results of our risk analysis in Section 4, and the corresponding proofs in Section 5, we will again make use of the function $\theta(\cdot)$ underlying $\boldsymbol{\theta}$ [although the role of $X(\cdot)$ will remain implicit through $\mathbf{X}$].

Let $\Theta$ be a generic function space to be defined later, such as the space of functions of bounded variation or a Besov space. We define our three models as follows.

(G) Gaussian model. Let $\theta \in \Theta$, and define $\theta_i = N \int_{I_i} \theta(t)\,dt$ to be the average of $\theta$ over $I_i$. Sample the $X_i$ independently as $X_i | \theta_i \sim \text{Gaussian}(\theta_i, \sigma^2)$, where $\sigma^2$ is assumed fixed and known.

The multiscale components in (4) then take the form

$$X_{\text{ch}(I),\text{l}} | X_I, \omega_I \sim \text{Gaussian}\left(\frac{N_{\text{ch}(I),\text{l}}}{N_I} X_I - \omega_I, c_I \sigma^2\right),$$

with

$$\omega_I = c_I \left(\frac{\theta_{\text{ch}(I),\text{r}}}{N_{\text{ch}(I),\text{r}}} - \frac{\theta_{\text{ch}(I),\text{l}}}{N_{\text{ch}(I),\text{l}}}\right),$$

for $c_I = N_{\text{ch}(I),\text{l}} N_{\text{ch}(I),\text{r}} / N_I$. The coarse scale component takes the form $X_{I_{00}} | \theta_{I_{00}} \sim \text{Gaussian}(\theta_{I_{00}}, N_{I_{00}} \sigma^2)$.

(P) Poisson model. Let $\theta \in \Theta$, where $\theta(t) \in [c, C]$, $\forall t \in [0, 1]$, for $0 < c < C$. Define $\theta_i = N \int_{I_i} \theta(t)\,dt$ to be the average of $\theta$ over $I_i$. Sample the $X_i$ independently as $X_i \sim \text{Poisson}(\theta_i)$.

The multiscale components in (4) then take the form

$$X_{\text{ch}(I),\text{l}} | X_I, \omega_I \sim \text{Binomial}(X_I; \omega_I),$$

with $\omega_I = \theta_{\text{ch}(I),\text{l}} / \theta_I$, while the coarse scale component takes the form $X_{I_{00}} | \theta_{I_{00}} \sim \text{Poisson}(\theta_{I_{00}})$.

(M) Multinomial model. Let $\theta \in \Theta$, where $\theta(t) \in [c, C]$ $\forall t \in [0, 1]$, for $0 < c < C$, and $\int_0^1 \theta(t)\,dt = 1$. Define $\theta_i = \int_{I_i} \theta(t)\,dt$. Let the components $X_i$ arise through (singular) multinomial sampling, that is, $\mathbf{X} \sim \text{Multinomial}(n; \boldsymbol{\theta})$, for some $n \sim N$.

The multiscale components in (4) then take the form

$$X_{\text{ch}(I),\text{l}} | X_I, \omega_I \sim \text{Binomial}(X_I; \omega_I),$$

with $\omega_I = \theta_{\text{ch}(I),\text{l}} / \theta_I$, as in the Poisson model, but the coarse scale component is now a point mass at $n$.

Model (G) is just the Gaussian "signal-plus-noise" model with average-sampling of the underlying function $\theta(\cdot)$, while model (P) is the Poisson analogue. Model (M) can be viewed as a discretized density estimation model, where a sample of size $n$ from the density $\theta$ is implicit. The criterion that $N$ be such that $n \sim N$ in this model may be satisfied, for example, by selecting $N$ to be the smallest power of 2 greater than or equal to $n$. This choice aids in producing computationally efficient implementations of the estimators defined in the next section, and has been found to work well in practice. The fact that the Poisson and multinomial models share the same multiscale components follows from the shared properties of the SSPSD class of distributional families, as explained above.

3.1. *Complexity penalized estimators.* The construction underlying the multiscale factorization in (4) involves intimate connections between factorizations, partitions and orthonormal bases, the exploitation of which is important for the creation of adaptive estimators and efficient algorithms for their calculation. That the factorization is closely linked to recursive partitioning is clear (i.e., through $\mathcal{P}^*$). However, through the latter, the former is also linked to a certain class of wavelet bases. To see this it is enough to note that, in the construction of any given C-RP $\mathcal{P}^*$, the splitting of each parent interval $I$ into its two children can be associated with a function $h_I$, as defined in (3). This function is a generalization of the dyadic Haar function $h_{j,k}$ defined earlier and the collection of all such $h_I$, for $I \in \mathcal{I}_{\mathrm{NT}}(\mathcal{P}^*)$, along with a single scale function on the full interval $[0, 1)$, are an example of what Girardi and Sweldens (1997) term an "unbalanced" Haar basis (UHB).

The multiscale coefficients $\omega_I$ in model (G) actually are proportional to the UHB coefficients $\langle \boldsymbol{\theta}, h_I \rangle$, differing only by their constants $c_I$ and $c_I'$. Therefore, in particular, $\omega_I = 0$ if and only if $\langle \boldsymbol{\theta}, h_I \rangle = 0$. On the other hand, in models (P) and (M) the $\omega_I$ arise as ratios of (left) child to parent sums of the appropriate components of $\boldsymbol{\theta}$. However, these ratios can be expressed as simple functions of the corresponding Haar coefficients, that is,

(6) $$\omega_I = \frac{\theta_{\mathrm{ch}(I),\mathrm{l}}}{\theta_I} = c_I' \bigg( \frac{c_I'}{N_{\mathrm{ch}(I),\mathrm{r}}} - \frac{\langle \boldsymbol{\theta}, h_I \rangle}{\theta_I} \bigg).$$

Note that $\omega_I = N_{\mathrm{ch}(I),\mathrm{l}}/N_I \equiv \rho_I$ if and only if $\langle \boldsymbol{\theta}, h_I \rangle = 0$. The value of $\rho_I$ is the ratio of the (discrete) lengths of the interval $I$ and its (left) child, and indicates homogeneity (smoothness) in $\boldsymbol{\theta}$ at the scale and position of $I$, just as $\omega_I = 0$ does in the Gaussian–wavelet case.

Hence the well-known exploitation of "sparseness" associated with wavelet expansions in Gaussian denoising problems also carries over to the Poisson and multinomial models, in that a piecewise-homogeneous vector $\boldsymbol{\theta}$ will have a large proportion of its $\omega_I$ equal to the corresponding $\rho_I$. This point suggests the promise of extensions of "keep-or-kill" thresholding algorithms to models (P) and (M), in which individual multiscale parameters $\omega_I$ are either set to their empirical value $X_{\mathrm{ch}(I),\mathrm{l}}/X_I$ (i.e., "keep") or to the default value $\rho_I$ (i.e., "kill"), as determined by whether a certain criterion function exceeds a threshold or not. In addition, in analogy to the results of Donoho (1997), in which the equivalence of a certain form of CART [Breiman, Friedman, Olshen and Stone (1984)] and algorithms based on constrained thresholding of appropriately defined Haar expansions are detailed, one might instead choose to base an estimator upon some optimal choice of recursive partition $\mathcal{P}$ among all such partitions $\mathcal{P} \preceq \mathcal{P}^*$, for some fixed choice of C-RP $\mathcal{P}^*$. For example, the complete recursive dyadic partition $\mathcal{P}^*_{\mathrm{Dy}}$ is a natural choice. Lastly, if one were to consider searching only with respect to a given C-RP $\mathcal{P}^*$ to be too constraining, a natural extension is to the entire library $\mathcal{L}$ of all $(N-1)!$ possible C-RP's $\mathcal{P}^*$ of the interval $[0, 1)$.

We therefore consider here estimators of thresholding, recursive dyadic partitioning and (general) recursive partitioning types, generically for each of our three models (G), (P) and (M). In the remaining sections we present results regarding the risk properties of certain simple versions of these types of estimators. Formally, we express the estimators in the form

$$(7) \qquad \hat{\boldsymbol{\theta}}(\mathbf{X}) \equiv \arg\min_{\boldsymbol{\theta}' \in \Gamma} \{-\log p(\mathbf{X}|\boldsymbol{\theta}') + 2\,\mathrm{pen}(\boldsymbol{\theta}')\},$$

$$(8) \qquad \mathrm{pen}(\boldsymbol{\theta}') \equiv \lambda \cdot \#\{\omega_I(\boldsymbol{\theta}') \text{ nontrivial}\},$$

where "nontrivial" means nonzero in the Gaussian case and not equal to $\rho_I$ in the Poisson and multinomial cases, $\lambda$ is a penalty factor to be defined later (e.g., Theorem 4) and $\Gamma$ is the space of possible values for a given estimator, defined as follows:

1. thresholding (T),

$$(9) \qquad \Gamma_{\mathrm{T}} \equiv \left\{ \boldsymbol{\theta}' \,\Big|\, \theta'_i = \beta_0 + \sum_{I \in \mathcal{I}} \beta_I h_I(i) \text{ for } \mathcal{I} \subseteq \mathcal{I}_{\mathrm{NT}}(\mathcal{P}^*_{\mathrm{Dy}}) \right\};$$

2. recursive dyadic partitioning (RDP),

$$(10) \qquad \Gamma_{\mathrm{RDP}} \equiv \left\{ \boldsymbol{\theta}' \,\Big|\, \theta'_i = \beta_0 + \sum_{I \in \mathcal{I}_{\mathrm{NT}}(\mathcal{P})} \beta_I h_I(i) \text{ for } \mathcal{P} \preceq \mathcal{P}^*_{\mathrm{Dy}} \right\};$$

3. recursive partitioning (RP),

$$(11) \qquad \Gamma_{\mathrm{RP}} \equiv \left\{ \boldsymbol{\theta}' \,\Big|\, \theta'_i = \beta_0 + \sum_{I \in \mathcal{I}_{\mathrm{NT}}(\mathcal{P})} \beta_I h_I(i) \text{ for } \mathcal{P} \preceq \mathcal{P}^*, \mathcal{P}^* \in \mathcal{L} \right\}.$$

In the definition of our spaces $\Gamma$ in (9)–(11), which we have expressed in terms of the UHB functions (3) for simplicity and comparison [and recall $\mathcal{I}_{\mathrm{NT}}(\mathcal{P})$ is the set of all nonterminal intervals encountered in the construction of a recursive partition $\mathcal{P}$], it is to be understood that the coefficient vectors $\boldsymbol{\beta}$ are constrained accordingly in each of the models (G), (P) and (M). That is, while there are no constraints in model (G), positivity of $\boldsymbol{\theta}'$ is required in model (P), and both positivity and unit summability in model (M). These latter two sets of constraints are enforced naturally in actual computations by virtue of working in the reparameterization $(\theta'_{I_{00}}, \omega'_I)$. Hence the multiscale reparameterizations are important not only algorithmically, through their role in the decoupling of terms in (4), but also mathematically in enforcing original constraints on $\boldsymbol{\theta}'$ in a natural manner.

Concerning the optimization in (7), by comparing the factorized form of the likelihood in (4) with the summability of the penalty in the same multiscale indexing in (8), it is not difficult to see that the estimator $\hat{\boldsymbol{\theta}}_\mathrm{T}$ is equivalent to performing a set of independent generalized likelihood ratio tests for each $I \in \mathcal{I}_\mathrm{NT}(\mathcal{P}^*_\mathrm{Dy})$. (Note that our choice of $\mathcal{P}^* = \mathcal{P}^*_\mathrm{Dy}$ here is arbitrary and made for convenience.) For model (G), since the log-likelihood is simply a sum of squares, this is in fact penalized least-squares with a counting penalty, which reduces to so-called hard-thresholding. On the other hand, for models (P) and (M) the index set $\mathcal{I}$ for the optimal estimate $\hat{\boldsymbol{\theta}}_\mathrm{T}$ corresponds to those indices $I$ for which the null hypothesis $H_I^{(0)}: \omega_I = \rho_I$ is rejected in favor of $H_I^{(1)}: \omega_I \neq \rho_I$, with respect to the local binomial likelihood functions. The estimator $\hat{\boldsymbol{\theta}}_\mathrm{RDP}$ is the analogue of the penalized least-squares estimator defined in Donoho (1997), wherein it is noted that recursive dyadic partitioning estimators like this are in fact thresholding estimators with hereditary constraints placed on which components may be "kept" or "killed" (i.e., resulting from the requirement that the partitioning be recursive). Finally, the estimator $\hat{\boldsymbol{\theta}}_\mathrm{RP}$ is an extension of this framework and reasoning to the larger space $\mathcal{L}$.

On a final note, we mention that all three of the estimators defined above may be computed in a computationally efficient manner, as summarized in the following.

THEOREM 3. *For models* (G), (P) *and* (M), *the following hold:*

(i) $\hat{\boldsymbol{\theta}}_\mathrm{T}$ *may be computed using an $O(N)$ thresholding procedure;*

(ii) $\hat{\boldsymbol{\theta}}_\mathrm{RDP}$ *may be computed using an $O(N)$ optimal tree pruning algorithm;*

(iii) $\hat{\boldsymbol{\theta}}_\mathrm{RP}$ *may be computed using an $O(N^3)$ dynamic programming algorithm.*

Proof of the theorem may be found in the Appendix.

**4. Risk analysis.** We state briefly in this section the main results deriving from a risk analysis for the estimators $\hat{\boldsymbol{\theta}}_T$, $\hat{\boldsymbol{\theta}}_{RDP}$ and $\hat{\boldsymbol{\theta}}_{RP}$ under the models (G), (P) and (M). In this section and the remainder of the paper, we will use $\mathbf{x}$ to denote an arbitrary element from the range of the random variable $\mathbf{X}$. We will measure the loss associated with estimating $\boldsymbol{\theta}$ by $\hat{\boldsymbol{\theta}}$ in terms of the (squared) Hellinger distance between the corresponding densities, that is,

$$
\begin{aligned}
L(\hat{\boldsymbol{\theta}}, \boldsymbol{\theta}) &\equiv H^2(p_{\hat{\boldsymbol{\theta}}}, p_{\boldsymbol{\theta}}) \\
&= \int \left[\sqrt{p(\mathbf{x}|\hat{\boldsymbol{\theta}})} - \sqrt{p(\mathbf{x}|\boldsymbol{\theta})}\right]^2 \nu(\mathbf{x}),
\end{aligned}
\tag{12}
$$

where $\nu$ is the dominating measure. The risk will be defined then as $R(\hat{\boldsymbol{\theta}}, \boldsymbol{\theta}) \equiv (1/N)E[L(\hat{\boldsymbol{\theta}}, \boldsymbol{\theta})]$, where the expectation is with respect to $\boldsymbol{\theta}$.

We assign properties to the values of the true $\boldsymbol{\theta}$ through properties of the function $\theta(\cdot)$ from which it was sampled. Given the role of Haar-like functions in our framework, a natural space $\Theta$ to consider is a ball $\Theta = \mathrm{BV}(C)$ of functions of bounded variation, that is, for which

$$
\sup_{M \geq 2} \sup_{t_1 \leq \cdots \leq t_M} \sum_{m=2}^{M} |\theta(t_m) - \theta(t_{m-1})| < C.
\tag{13}
$$

We then have the following result regarding upper bounds on the risk.

THEOREM 4. *Assume $\Theta = \mathrm{BV}(C)$ and the conditions of either model (G), (P) or (M). Let the constant $\lambda$ in (8) be of the form $\gamma \log N$, for $\gamma \geq 3/2$ and $N \geq 3$. Then the risks of the estimators $\hat{\boldsymbol{\theta}}_T$ and $\hat{\boldsymbol{\theta}}_{RP}$ are bounded above by $O((\log N/N)^{2/3})$, while the risk of the estimator $\hat{\boldsymbol{\theta}}_{RDP}$ is bounded above by $O((\log^2 N/N)^{2/3})$.*

Proof of this result and all others given in this section can be found in Section 5. Note that the performances of the estimators are bounded identically for the three models. A minor variation in the overall proof leads to risk statements similar to those of Theorem 4 when loss functions of squared-error type are used.

COROLLARY 1. *The same upper bounds hold when risk is measured for the models (G), (P) and (M) as $(1/4N)E[\|\boldsymbol{\theta} - \hat{\boldsymbol{\theta}}\|^2]$, $(1/N)E[\|\boldsymbol{\theta}^{1/2} - \hat{\boldsymbol{\theta}}^{1/2}\|^2]$ and $(1/N)E[\|(n\boldsymbol{\theta})^{1/2} - (n\hat{\boldsymbol{\theta}})^{1/2}\|^2]$, respectively.*

As mentioned in the Introduction, two key properties of wavelet-based thresholding estimators are their near-optimal risk behavior and their adaptivity. To establish near-optimality of our multiscale estimators we need the following lower bound on the minimax risk.

THEOREM 5. *Assume the conditions of Theorem 4 and Hellinger loss. Then the minimax risk obeys a lower bound of $O(N^{-2/3})$ in each of the models* (G), (P) *and* (M).

Hence, in particular, combining the upper and lower bounds in Theorems 4 and 5, we obtain that the estimators $\hat{\boldsymbol{\theta}}_{\mathrm{T}}$ and $\hat{\boldsymbol{\theta}}_{\mathrm{RP}}$ come within the same logarithmic factor of minimax risk, while the estimator $\hat{\boldsymbol{\theta}}_{\mathrm{RDP}}$ comes within the square of that factor, for each of the models (G), (P) and (M).

Adaptivity of our estimators then follows from the fact that similar near-optimality statements hold in other spaces of varying smoothness, despite the fact that the estimators have no a priori knowledge of which space(s) contains $\theta(\cdot)$. For example, for an appropriately defined range of Besov spaces $B_{p,q}^{\xi}$ we have the following.

THEOREM 6. *Suppose $\Theta = B_{p,q}^{\xi}([0,1])$ is a Besov space, with $0 < \xi < 1$ and $1 \leq p < \infty$ such that $1/p < \xi + 1/2$, and $q > 0$. Then the conclusions of Theorems 4 and 5, as well as Corollary 1, hold with the exponent $2/3$ replaced by $2\xi/(2\xi + 1)$.*

In summary, the results of this section describe, for smoothness classes appropriate to Haar-like bases, how all three of our multiscale complexity penalized likelihood estimators exhibit the same sort of adaptivity and near-optimality properties as classical wavelet-based models—but simultaneously for certain continuous, count and categorical data types. As the proofs in the next section demonstrate, this success is due primarily to (i) the role of Haar-like multiresolution structures in our framework and (ii) the ability to decouple the information in the data through likelihood factorizations that mirror these structures. Extensions of these results to classes of smoother objects (e.g., Besov spaces with $\xi > 1$) would seem to require wavelets smoother than those in a Haar basis. However, achieving similar likelihood factorizations in this context is a much more difficult problem, and one for which it seems unrealistic to expect the type of complete multiscale decoupling of information in data and parameters inherent in (4).

**5. Proof of main results.** We establish proofs of the results in Section 4 here. Although the case of model (G) could be handled using arguments from the existing literature on wavelets and Gaussian noise models, such arguments do not immediately extend to the cases of models (P) and (M). Hence, since our interest is in part to illustrate how all three models may be handled in a unified fashion, we introduce a generally applicable result in Section 5.1. The results for model (G), within this general framework, are then presented in Section 5.2, while the results for models (P) and (M) follow in Section 5.3.

5.1. *A fundamental risk bound.* Our proof of Theorem 4, and hence the near-optimality of our estimators, rests primarily upon a fundamental upper bound on the expected (squared) Hellinger distance. The form of this bound is much like those in Barron, Birgé and Massart (1999). However, whereas the proof of their (quite general) bounds relies on recent advances in isoperimetric inequalities and a great deal of careful technical work, our own bounds adapt recent arguments of Li (1999) and Li and Barron (2000) for mixture-based density estimation which, in particular, rely on the discretization (quantization) of the space $\Gamma$ of estimates [in the spirit of earlier work by Barron and Cover (1991)] to produce a relatively simple proof, applicable to all three models under consideration here.

Let $H^2(p_{\boldsymbol{\theta}^{(1)}}, p_{\boldsymbol{\theta}^{(2)}})$ be the (squared) Hellinger distance between two densities $p(\mathbf{x}|\boldsymbol{\theta}^{(1)})$ and $p(\mathbf{x}|\boldsymbol{\theta}^{(2)})$, as introduced earlier. Additionally, define the Kullback–Leibler divergence between these densities as

$$(14) \qquad K(p_{\boldsymbol{\theta}^{(1)}}, p_{\boldsymbol{\theta}^{(2)}}) = \int \log \frac{p(\mathbf{x}|\boldsymbol{\theta}^{(1)})}{p(\mathbf{x}|\boldsymbol{\theta}^{(2)})} p(\mathbf{x}|\boldsymbol{\theta}^{(1)}) \nu(\mathbf{x}).$$

The following theorem bounds the expected (squared) Hellinger distance in terms of the Kullback–Leibler divergence.

THEOREM 7. *Let $\Gamma_N$ be a finite collection of estimators $\boldsymbol{\theta}'$ for $\boldsymbol{\theta}$, and pen$(\cdot)$ a function on $\Gamma_N$ satisfying the condition*

$$(15) \qquad \sum_{\boldsymbol{\theta}' \in \Gamma_N} e^{-\operatorname{pen}(\boldsymbol{\theta}')} \leq 1.$$

*Let $\hat{\boldsymbol{\theta}}$ be a penalized maximum likelihood estimator of the form*

$$(16) \qquad \hat{\boldsymbol{\theta}}(\mathbf{X}) \equiv \arg \min_{\boldsymbol{\theta}' \in \Gamma_N} \{-\log p(\mathbf{X}|\boldsymbol{\theta}') + 2\operatorname{pen}(\boldsymbol{\theta}')\}.$$

*Then*

$$(17) \qquad E[H^2(p_{\hat{\boldsymbol{\theta}}}, p_{\boldsymbol{\theta}})] \leq \min_{\boldsymbol{\theta}' \in \Gamma_N} \{K(p_{\boldsymbol{\theta}}, p_{\boldsymbol{\theta}'}) + 2\operatorname{pen}(\boldsymbol{\theta}')\}.$$

PROOF. The proof uses the same key ideas as found in Li (1999) and Li and Barron (2000). Begin by noting that

$$(18) \qquad \begin{aligned} H^2(p_{\boldsymbol{\theta}^{(1)}}, p_{\boldsymbol{\theta}^{(2)}}) &= \int \left[\sqrt{p(\mathbf{x}|\boldsymbol{\theta}^{(1)})} - \sqrt{p(\mathbf{x}|\boldsymbol{\theta}^{(2)})}\right]^2 \nu(\mathbf{x}) \\ &= 2\left(1 - \int \sqrt{p(\mathbf{x}|\boldsymbol{\theta}^{(1)}) p(\mathbf{x}|\boldsymbol{\theta}^{(2)})} \nu(\mathbf{x})\right) \\ &\leq -2\log \int \sqrt{p(\mathbf{x}|\boldsymbol{\theta}^{(1)}) p(\mathbf{x}|\boldsymbol{\theta}^{(2)})} \nu(\mathbf{x}), \end{aligned}$$

where $\nu$ is the dominating measure. Taking the expectation with respect to $\mathbf{X}$, we then have

$$E[H^2(p_{\hat{\boldsymbol{\theta}}}, p_{\boldsymbol{\theta}})] \leq 2E\log\left(\frac{1}{\int[\sqrt{p(\mathbf{x}|\hat{\boldsymbol{\theta}})p(\mathbf{x}|\boldsymbol{\theta})}]\,\nu(\mathbf{x})}\right)$$

$$\leq 2E\log\left(\frac{p^{1/2}(\mathbf{X}|\hat{\boldsymbol{\theta}})e^{-\operatorname{pen}(\hat{\boldsymbol{\theta}})}}{p^{1/2}(\mathbf{X}|\tilde{\boldsymbol{\theta}})e^{-\operatorname{pen}(\tilde{\boldsymbol{\theta}})}}\frac{1}{\int[\sqrt{p(\mathbf{x}|\hat{\boldsymbol{\theta}})p(\mathbf{x}|\boldsymbol{\theta})}]\,\nu(\mathbf{x})}\right),$$

where $\tilde{\boldsymbol{\theta}}$ is the argument that minimizes the right-hand side of the expression in (17), which is simply the theoretical analogue of $\hat{\boldsymbol{\theta}}$. However, the last expression above may be broken into two pieces, being equal to

(19) $$E\left[\log\frac{p(\mathbf{X}|\boldsymbol{\theta})}{p(\mathbf{X}|\tilde{\boldsymbol{\theta}})}\right] + 2\operatorname{pen}(\tilde{\boldsymbol{\theta}})$$

(20) $$+ 2E\log\left(\frac{p^{1/2}(\mathbf{X}|\hat{\boldsymbol{\theta}})}{p^{1/2}(\mathbf{X}|\boldsymbol{\theta})}\frac{e^{-\operatorname{pen}(\hat{\boldsymbol{\theta}})}}{\int[\sqrt{p(\mathbf{x}|\hat{\boldsymbol{\theta}})p(\mathbf{x}|\boldsymbol{\theta})}]\,\nu(\mathbf{x})}\right).$$

Noting that the expression in (19) is just the right-hand side of (17), the rest of the proof entails showing that the expression in (20) is bounded above by zero. Specifically, applying Jensen's inequality we have the upper bound

(21) $$2\log E\left[e^{-\operatorname{pen}(\hat{\boldsymbol{\theta}})}\frac{\sqrt{(p(\mathbf{X}|\hat{\boldsymbol{\theta}})/p(\mathbf{X}|\boldsymbol{\theta}))}}{\int[\sqrt{p(\mathbf{x}|\hat{\boldsymbol{\theta}})p(\mathbf{x}|\boldsymbol{\theta})}]\,\nu(\mathbf{x})}\right].$$

The integrand in the expectation in (21) can be upper bounded by

(22) $$\sum_{\boldsymbol{\theta}'\in\Gamma_N} e^{-\operatorname{pen}(\boldsymbol{\theta}')}\frac{\sqrt{(p(\mathbf{X}|\boldsymbol{\theta}')/p(\mathbf{X}|\boldsymbol{\theta}))}}{\int[\sqrt{p(\mathbf{x}|\boldsymbol{\theta}')p(\mathbf{x}|\boldsymbol{\theta})}]\,\nu(\mathbf{x})},$$

which, since $\boldsymbol{\theta}'$ does not depend on $\mathbf{X}$, produces the following upper bound for (20):

(23) $$2\log\sum_{\boldsymbol{\theta}'\in\Gamma_N} e^{-\operatorname{pen}(\boldsymbol{\theta}')}\frac{E[\sqrt{(p(\mathbf{X}|\boldsymbol{\theta}')/p(\mathbf{X}|\boldsymbol{\theta}))}\,]}{\int[\sqrt{p(\mathbf{x}|\boldsymbol{\theta}')p(\mathbf{x}|\boldsymbol{\theta})}]\,\nu(\mathbf{x})}.$$

Now the integration $\int[\sqrt{p(\mathbf{x}|\boldsymbol{\theta}')p(\mathbf{x}|\boldsymbol{\theta})}]\,\nu(\mathbf{x})$ in the denominator can be rewritten as an expectation with respect to the distribution of $\mathbf{X}$ through multiplication by $p(\mathbf{x}|\boldsymbol{\theta})/p(\mathbf{x}|\boldsymbol{\theta})$. However, this yields the same expectation as appears in the numerator of each term in (23), and hence (23) is equal to $2\log\sum_{\boldsymbol{\theta}'\in\Gamma_N} e^{-\operatorname{pen}(\boldsymbol{\theta}')}$. This and the condition in (15) yield that (23) is bounded by zero. $\square$

The result of Theorem 7 holds quite generally, and in particular for each of the densities of models (G), (P) and (M). Use of this bound for statements regarding the estimators $\hat{\boldsymbol{\theta}}_{\mathrm{T}}$, $\hat{\boldsymbol{\theta}}_{\mathrm{RDP}}$ and $\hat{\boldsymbol{\theta}}_{\mathrm{RP}}$ requires that the inequality in (15) holds, for which it is sufficient to establish the following.

LEMMA 1. *Let $\Gamma_N$ be the collection of all $N$-length vectors $\boldsymbol{\theta}'$ with components $\theta'_i \in D_N[R_1, R_2]$, for some $R_1 < R_2$, where $D_N[R_1, R_2]$ denotes a discretization of the interval $[R_1, R_2]$ into $N^{1/2}$ equispaced values. Let $\#(\boldsymbol{\theta}')$ count the number of constant-valued sequences in the vector $\boldsymbol{\theta}'$, that is, in analogy to the number of pieces of a piecewise constant function. Then*

$$\sum_{\boldsymbol{\theta}' \in \Gamma_N} e^{-\gamma \log N \#(\boldsymbol{\theta}')} \leq 1, \tag{24}$$

*for $\gamma \geq 3/2$ and $N \geq 3$.*

PROOF. Begin by writing $\Gamma_N = \bigcup_{d=1}^{N} \Gamma_N^{(d)}$, where $\Gamma_N^{(d)}$ is the subset of values $\boldsymbol{\theta}'$ that is composed of $d$ constant-valued sequences. For example, $(1, 1, 2, 2, 3)$ and $(1, 2, 3, 3, 3)$ might be two such sequences in $\Gamma_5^{(3)}$. Each of the members of $\Gamma_N^{(d)}$ has the same value for the summand in (24), and there are $(N^{1/2})^d$ distinct values that may be taken on by the set of $d$ constant-valued components of each member. Also, there are $N - 1$ choose $d - 1$ possible $d$-component vectors of length $N$. So we have

$$\begin{aligned}
\sum_{\boldsymbol{\theta}' \in \Gamma_N} e^{-\gamma \log N \#(\boldsymbol{\theta}')} &= \sum_{d=1}^{N} \binom{N-1}{d-1} e^{-(\gamma - 1/2) d \log N} \\
&= \sum_{d'=0}^{N-1} \binom{N-1}{d'} e^{-(\gamma - 1/2)(d'+1) \log N} \\
&\leq \sum_{d'=0}^{N-1} \frac{(N-1)^{d'}}{d'!} N^{-(\gamma - 1/2)(d'+1)} \\
&\leq N^{-(\gamma - 1/2)} \sum_{d'=0}^{N-1} \frac{1}{d'!} \leq N^{-(\gamma - 1/2)} e,
\end{aligned}$$

which is bounded by 1 under the conditions given. □

In the cases of $\hat{\boldsymbol{\theta}}_{\mathrm{T}}$ and $\hat{\boldsymbol{\theta}}_{\mathrm{RP}}$, the number of nontrivial $\omega_I$ defining the penalty in (8) is in fact the same as the penalty in the statement of the lemma. The spaces $\Gamma_{\mathrm{T}}$ and $\Gamma_{\mathrm{RP}}$, with appropriate discretization (to be defined below), are equivalent to $\Gamma_N$. For the case of $\hat{\boldsymbol{\theta}}_{\mathrm{RDP}}$, it is enough to note the inclusion $\Gamma_{\mathrm{RDP}} \subset \Gamma_{\mathrm{RP}}$.

### 5.2. *The Gaussian case.*

PROOF OF THEOREM 4—MODEL (G). Assume without loss of generality that $\sigma \equiv 1$. Using inequality (17) and the fact that the Kullback–Leibler divergence in the Gaussian case is simply proportional to a (squared) $\ell_2$-norm, we have that

$$(25) \qquad R(\hat{\boldsymbol{\theta}}, \boldsymbol{\theta}) \leq \min_{\boldsymbol{\theta}' \in \Gamma_N} \left\{ \frac{1}{2N} \|\boldsymbol{\theta} - \boldsymbol{\theta}'\|_{\ell_2}^2 + \frac{2\gamma \log N}{N} \#(\boldsymbol{\theta}') \right\}.$$

The minimization in (25) essentially seeks an optimal balancing of bias and variance terms, respectively. We will bound this quantity by bounding the bias term over $\Gamma_N^{(d)}$, for each fixed $d$, and then optimizing our resulting overall bound in $d$. In producing a bound on the bias, the following result from Donoho (1993) is central.

LEMMA 2. *Let $\theta(\cdot) \in \mathrm{BV}$. Define $\tilde{\theta}(\cdot)$ to be the best $d$-term approximant to $\theta(\cdot)$ in the dyadic Haar basis for $L_2([0,1])$. Then $\|\theta - \tilde{\theta}\|_{L_2} = O(1/d)$.*

Define $\tilde{\boldsymbol{\theta}}$ to be the average sampling of $\tilde{\theta}$ on the intervals $I_i$, where the dependence of $\tilde{\boldsymbol{\theta}}$ on $d$ is to be understood. Then let $\tilde{\boldsymbol{\theta}}' \equiv [\tilde{\boldsymbol{\theta}}]$ be the result of quantizing the elements of $\tilde{\boldsymbol{\theta}}$ to the set $D_N[-C, C]$, where $C$ is the radius of the BV ball in the statement of Theorem 4. By the triangle inequality it follows that

$$(26) \quad \begin{aligned} (1/N)\|\boldsymbol{\theta} - \tilde{\boldsymbol{\theta}}'\|_{\ell_2}^2 &\leq (1/N)\|\boldsymbol{\theta} - \tilde{\boldsymbol{\theta}}\|_{\ell_2}^2 + (1/N)\|\tilde{\boldsymbol{\theta}} - \tilde{\boldsymbol{\theta}}'\|_{\ell_2}^2 \\ &\quad + (2/N)\|\boldsymbol{\theta} - \tilde{\boldsymbol{\theta}}\|_{\ell_2} \|\tilde{\boldsymbol{\theta}} - \tilde{\boldsymbol{\theta}}'\|_{\ell_2}. \end{aligned}$$

The first term on the right-hand side is a measure of approximation error in $\ell_2(N)$, which may be bounded by the corresponding $L_2$ approximation error by exploiting the piecewise constant nature of the Haar basis and our use of average sampling. Specifically, let $h_{j,k}(i)$ be the $(j,k)$th Haar function on the discrete space $\{0, 1, \ldots, N-1\}$, as defined in Section 2, and let $h_{j,k}^c(t)$ be its analogue on the interval $[0, 1]$. Then it follows that $\langle \boldsymbol{\theta}, h_{j,k} \rangle_{\ell_2} = N^{1/2} \langle \theta, h_{j,k}^c \rangle_{L_2}$ and therefore

$$(27) \quad \begin{aligned} (1/N)\|\boldsymbol{\theta} - \tilde{\boldsymbol{\theta}}\|_{\ell_2}^2 &= (1/N) \sum_{(j,k) \in \mathcal{J}_N} (\langle \boldsymbol{\theta}, h_{j,k} \rangle_{\ell_2} - \langle \tilde{\boldsymbol{\theta}}, h_{j,k} \rangle_{\ell_2})^2 \\ &= \sum_{(j,k) \in \mathcal{J}_N} (\langle \theta, h_{j,k}^c \rangle_{L_2} - \langle \tilde{\theta}, h_{j,k}^c \rangle_{L_2})^2, \end{aligned}$$

where $\mathcal{J}_N$ is the set of $(j, k)$ with $j = 0, 1, \ldots, J-1$ and $k = 0, 1, \ldots, 2^j - 1$. However, the last term in (27) is bounded by a similar sum over all $(j, k)$, which in turn is equal to the (squared) $L_2$ approximation error $\|\theta - \tilde{\theta}\|_{L_2}^2$.

Hence the first term in (26) is of order $O(d^{-2})$. The second term is simply a discretization error, which is controlled through the conditions of Lemma 1 to be of order $O(1/N)$. Therefore the cross-term in (26) is of order $O(d^{-1}N^{-1/2})$. In the case of estimation through the thresholding or recursive partitioning strategies, the quantity $\#(\widetilde{\boldsymbol{\theta}})$ will be proportional to $d$. Combining the above results [and ignoring the negligible $O(1/N)$ term] yields the bound

$$(28) \quad \min_{\boldsymbol{\theta}' \in \Gamma_N^{(d)}} \left\{ \frac{1}{2N} \|\boldsymbol{\theta} - \boldsymbol{\theta}'\|_{\ell_2}^2 + \frac{2\gamma \log N}{N} \#(\boldsymbol{\theta}') \right\} \\ \leq O(d^{-2}) + O(d^{-1}N^{-1/2}) + O(dN^{-1}\log N),$$

which is minimized for $d \sim (N/\log N)^{1/3}$. Substitution then yields the result that $R(\hat{\boldsymbol{\theta}}, \boldsymbol{\theta})$ is bounded by a quantity of order $O((\log N/N)^{2/3})$. A similar argument holds for estimation via recursive dyadic partitioning (RDP), where $\#(\widetilde{\boldsymbol{\theta}})$ instead behaves like $d \log N$, yielding the bound $O((\log_e^2(N)/N)^{2/3})$. □

PROOF OF COROLLARY 1—MODEL (G). Proof of this corollary, for all three models, is inherent in the proof of Theorem 7. Specifically, following Li (1999), define the "affinity" between two densities $p$ and $q$ as $\mathcal{A}(p,q) \equiv \int (p(\mathbf{x})q(\mathbf{x}))^{1/2} \nu(\mathbf{x})$. Then (18) can be rewritten as

$$(29) \quad H^2(p_{\boldsymbol{\theta}^{(1)}}, p_{\boldsymbol{\theta}^{(2)}}) \leq -2 \log \mathcal{A}(p_{\boldsymbol{\theta}^{(1)}}, p_{\boldsymbol{\theta}^{(2)}}),$$

and therefore Theorem 7 equivalently can be viewed as a bound on minus twice the log-affinity and related quantities thereof. For example, under independent sampling we have that $\mathcal{A}(p_{\boldsymbol{\theta}^{(1)}}, p_{\boldsymbol{\theta}^{(2)}}) = \prod_i \mathcal{A}(p_{\theta_i^{(1)}}, p_{\theta_i^{(2)}})$, and a short calculation shows that for model (G) $-2 \log \mathcal{A}(p_{\theta_i}, p_{\hat{\theta}_i}) = (1/4)(\theta_i - \hat{\theta}_i)^2$. □

PROOF OF THEOREM 5—MODEL (G). In the Gaussian case this follows from standard arguments, as have been used for analogous statements for wavelet-based estimators with the Gaussian signal-plus-noise model. See Donoho (1993), for example. The approach is based on the so-called method of hyperrectangles of Donoho, Liu and MacGibbon (1990). That is, one defines an object

$$(30) \quad \mathcal{H}_j \equiv \left\{ \sum_{k=0}^{2^j - 1} \beta_{j,k} h_{j,k}^c(t) : |\beta_k| \leq \Delta_j \right\},$$

where the choice of $\Delta_j \propto 2^{-3j/2}$ is made to produce a hypercube "just barely" in our BV ball $\Theta = \text{BV}(C)$, and $j = j^*$ is chosen to satisfy the constraint $2^{3j/2} \sim N^{1/2}$ so as to produce the hardest possible estimation problem on such hyperrectangles. Because the $\ell_2$ risk for estimating objects in $\mathcal{H}_j$ is simply the sum of the $\ell_2$ risks for estimating the individual $\beta_{j,k}$ in this setting, and because this latter can be bounded below by $2^{j^*}\varepsilon^2$, where $\varepsilon \propto N^{-1/2}$ is the noise level, the lower bound of $N^{-2/3}$ on the minimax risk follows. □

PROOF OF THEOREM 6—MODEL (G). The analogue of Theorem 4 in this context is proven simply by modifying the statement of Lemma 2 accordingly. That is, for Besov spaces in the range specified, the best $d$-term approximant in the dyadic Haar basis for $L_2([0,1])$ has an error $\|\theta - \tilde{\theta}\|_{L_2}$ of order $O(d^{-\xi})$. See DeVore (1998), for example. The effect of this change is to change the lead term in (28) from $O(d^{-2})$ to $O(d^{-2\xi})$, and the cross-term similarly, from which the upper bounds follow. Corollary 1 then follows as before. Finally, proof of Theorem 5 follows as in the case of BV, but with the appropriate changes made to the definition of $\Delta$ and $j^*$. □

### 5.3. *The Poisson and multinomial cases.*

PROOF OF THEOREM 4—MODELS (P) AND (M). As remarked previously, the bound in Theorem 7 holds for these two models as well. However, the Kullback–Leibler divergence in these cases is not equivalent to an $\ell_2$-norm. Nevertheless, we may pursue a strategy similar to that in the Gaussian case by bounding the size of the Kullback–Leibler term when evaluated at a particular estimate relating to an optimal nonlinear Haar approximant associated with the function $\theta(\cdot)$.

Begin with model (P). Let $\tilde{\theta}(\cdot)$ be the best $d$-term nonlinear approximant to $\theta(\cdot)$, in the sense of Lemma 2, and define $\tilde{\theta}_i = N \int_{I_i} \tilde{\theta}(t)\,dt$ through average-sampling. Now the condition that the function $\theta(\cdot) \in [c, C]$ and the use of average-sampling ensure that for the elements of $\boldsymbol{\theta}$ we have $\theta_i \in [c, C]$ as well. However, we also have that $\tilde{\theta}_i \in [c, C]$. This results from the fact that the $\tilde{\theta}_i$ derive from average-sampling the function $\tilde{\theta}(\cdot)$, and this latter has a range restricted to $[c, C]$. This last statement follows from noting that, by definition of minimizing $\|\theta - f\|_{L_2(0,1)}$ in a Haar basis, the function $\tilde{\theta}$ equivalently is defined by a set of characteristic functions on some dyadic subintervals $I$ that partition $[0, 1]$ and associated constants $\alpha_I = |I|^{-1} \int_I \theta(t)\,dt$. The $\alpha_I$ clearly are bounded by $c$ and $C$.

Continuing, let $\tilde{\boldsymbol{\theta}}' \equiv [\tilde{\boldsymbol{\theta}}]$ be the result of quantizing the elements of $\tilde{\boldsymbol{\theta}}$ to the set $D_N[c, C]$. Then the Kullback–Leibler divergence may be bounded as

$$\frac{1}{N} K(p_{\boldsymbol{\theta}}, p_{\tilde{\boldsymbol{\theta}}'}) = \frac{1}{N} \sum_{i=0}^{N-1} \tilde{\theta}'_i - \theta_i + \theta_i \log\left(\frac{\theta_i}{\tilde{\theta}'_i}\right)$$

$$\leq \frac{1}{N} \sum_{i=0}^{N-1} \tilde{\theta}'_i - \theta_i + \theta_i \left( \frac{\theta_i}{\tilde{\theta}'_i} - 1 \right)$$

(31)
$$= \frac{1}{N} \sum_{i=0}^{N-1} \left( \frac{1}{\tilde{\theta}'_i} \right) (\tilde{\theta}'^2_i + \theta_i^2 - 2\tilde{\theta}'_i \theta_i)$$

$$\leq \frac{1}{Nc} \|\boldsymbol{\theta} - \tilde{\boldsymbol{\theta}}'\|_{\ell_2}^2,$$

where the first inequality follows from $\log(z) \leq z - 1$ and the second follows from the fact that $\tilde{\theta}'_i \in [c, C] \ \forall i$. The inequality in (31) leaves us in essentially the same position with which we started in the Gaussian case. Hence, arguing in an entirely analogous fashion we are led to the same inequality as in (28), and the result of Theorem 4 is proven for model (P).

The argument for model (M) is similar in structure to that for model (P), but differs with respect to certain important technical details, deriving from the fact that $\tilde{\boldsymbol{\theta}}'$ must be both positive and sum to 1. Specifically, we begin by writing $\theta(t) = 1 + (\theta(t) - 1) \equiv 1 + g(t)$. Then, let $\tilde{g}(\cdot)$ be the best $d$-piece nonlinear Haar approximant to $g(\cdot)$. Since $\theta \in \Theta \Rightarrow g \in \Theta$, it follows that $\|g - \tilde{g}\|_{L_2} = O(d^{-1})$.

Next note that the Haar scale coefficient of $g(\cdot)$ at $j = 0$ is zero, from which it follows that that scale coefficient will be zero for $\tilde{g}(\cdot)$ as well, and therefore the latter is defined purely in terms of $d$ nonzero wavelet coefficients. Since the wavelets have zero integral, defining $\tilde{\theta}(t) \equiv 1 + \tilde{g}(t)$ results in an approximant to $\theta$ that integrates to 1. Hence $\tilde{\theta}(\cdot)$ will be a proper density if it is nonnegative. However, an argument similar to that of the Poisson case can be used to show that $g \in (c-1, C-1)$ implies the same for $\tilde{g}(\cdot)$, from which it follows that $\tilde{\theta} \in (c, C)$.

Define $\tilde{\boldsymbol{\theta}}$ through integration (i.e., *not* average-sampling) via $\tilde{\theta}_i \equiv \int_{I_i} \tilde{\theta}(t)\, dt$ and note that $\tilde{\boldsymbol{\theta}}$ will be a proper probability mass function on the set $\{0, 1, \ldots, N-1\}$, with elements $\tilde{\theta}_i$ bounded below by $c/N$. It remains for us to quantize $\tilde{\boldsymbol{\theta}}$ in such a manner as to preserve this property, which we accomplish by working instead with $N\tilde{\boldsymbol{\theta}}$. Noting that $N\tilde{\theta}_i \in [c, C] \ \forall i$, we quantize each of these elements away from zero in the positive (i.e., toward $+\infty$) direction on $D_N[c, C]$. Similar to the Gaussian and Poisson cases, this means that $|[N\tilde{\theta}_i] - N\tilde{\theta}_i| \sim N^{-1/2}$, by definition of $D_N[c, C]$.

Division of $[N\tilde{\boldsymbol{\theta}}]$ by $N$ would produce an object on the same scale as $\tilde{\boldsymbol{\theta}}$, with components in $[c, C]$, but it would no longer be a proper probability mass function because our method of quantization leads to an inflation of the overall mass by the amount $\mu \equiv -N + \sum_{i=0}^{N-1} [N\tilde{\theta}_i]$. We can correct for this increase by subtracting $\delta \equiv \mu/N \sim O(N^{-1/2})$ from each element of $[N\tilde{\boldsymbol{\theta}}]$. For $N$ sufficiently large, say $c/2 > N^{-1/2}$, our final estimator $\tilde{\boldsymbol{\theta}}' \equiv ([N\tilde{\boldsymbol{\theta}}] - \delta)/N$ is a proper probability mass function and is bounded below by $c/(2N)$.

Now, similar to (31), bound the Kullback–Leibler divergence between $p_{\boldsymbol{\theta}}$ and $p_{\tilde{\boldsymbol{\theta}}'}$ as

$$\frac{1}{N}K(p_{\boldsymbol{\theta}}, p_{\tilde{\boldsymbol{\theta}}'}) = \frac{1}{N}\sum_{i=0}^{N-1} n\theta_i \log\left(\frac{\theta_i}{\tilde{\theta}'_i}\right)$$

(32)
$$= \frac{n}{N}\sum_{i=0}^{N-1} \tilde{\theta}'_i - \theta_i + \theta_i \log\left(\frac{\theta_i}{\tilde{\theta}'_i}\right)$$

$$\leq \frac{n}{N}\frac{2N}{c}\|\boldsymbol{\theta} - \tilde{\boldsymbol{\theta}}'\|_{\ell_2}^2,$$

where the factor of $N$ in the numerator of the last line comes from $\tilde{\theta}'_i \geq \frac{c}{2N}$. Under the condition $n \sim N$ in model (M), we are left with something that behaves like $N\|\boldsymbol{\theta} - \tilde{\boldsymbol{\theta}}'\|_{\ell_2}^2$ in (32). Similar to the previous cases, this may be bounded by terms involving approximation error, discretization (quantization) error and a cross term, using the triangle inequality. Noting that the use of simple integration (as opposed to average-sampling) here leads to the alternative relation $\langle \boldsymbol{\theta}, h_{j,k}\rangle_{\ell_2} = N^{-1/2}\langle \theta, h_{j,k}^c\rangle_{L_2}$ between the discrete and continuous Haar coefficients, we find that the quantity $N\|\boldsymbol{\theta} - \tilde{\boldsymbol{\theta}}\|_{\ell_2}^2$ may be bounded above by $\|\theta - \tilde{\theta}\|_{L_2}^2 = \|g - \tilde{g}\|_{L_2}^2 = O(d^{-2})$. Similarly, by construction we have $\|\tilde{\boldsymbol{\theta}} - \tilde{\boldsymbol{\theta}}'\|_{\ell_2}^2 = O(N^{-2})$, and so the discretization error is again $O(N^{-1})$. Therefore, an inequality like that in (28) holds, and the result of Theorem 4 is proven for model (M). □

Note that, as a consequence of our argument, the condition $n \sim N$ arises in a natural manner. The interpretation of this condition is, viewed from the context of density estimation, that the number of total possible bins $N$ should be chosen on the order of the number of samples $n$. The underlying algorithms will choose an optimal number less than or equal to $n$ (in fact, likely much less), due to the fact that the case $N > n$ assures that there will be empty bins and these will be aggregated over.

PROOF OF COROLLARY 1—MODELS (P) AND (M). Model (P) involves independent sampling, and thus it suffices to note that a short calculation produces the expression $-2\log \mathcal{A}(p_{\theta_i}, p_{\hat{\theta}_i}) = (\theta_i^{1/2} - \hat{\theta}_i^{1/2})^2$. Now consider model (M), where we have

$$\mathcal{A}(p_{\boldsymbol{\theta}}, p_{\hat{\boldsymbol{\theta}}}) = \sum_{\mathbf{x}:\sum x_i = n} \binom{n}{x_0, \ldots, x_{N-1}} \prod_{i=0}^{N-1} (\theta_i \hat{\theta}_i)^{x_i/2} = \left(\sum_{i=0}^{N-1} (\theta_i \hat{\theta}_i)^{1/2}\right)^n.$$

The right-hand side of the above equation is itself an affinity to the power $n$, since both $\boldsymbol{\theta}$ and $\hat{\boldsymbol{\theta}}$ sum to 1. Therefore minus twice the logarithm of this quantity is an upper bound on $n$ times the (squared) Hellinger distance between these two vectors. $\square$

PROOF OF THEOREM 5—MODELS (P) AND (M). Our argument here is in the spirit of the method of orthogonal hyperrectangles outlined in the Gaussian case, but with a number of technical differences. Consider the Poisson case first and begin by noting that it is not difficult to show that the constraint $\theta(\cdot), \theta'(\cdot) \in [c, C]$ implies that the Kullback–Leibler divergence and (squared) Hellinger distance between densities corresponding to $\boldsymbol{\theta}$ and $\boldsymbol{\theta}'$ are within a constant factor of each other, where the constant is a function of $c$ and $C$. Hence it suffices to provide a lower bound on the quantity

$$\inf_{\hat{\boldsymbol{\theta}}} \sup_{\boldsymbol{\theta}} \frac{1}{N} E[K(p_{\boldsymbol{\theta}}, p_{\hat{\boldsymbol{\theta}}})]. \tag{33}$$

Because the Kullback–Leibler divergence is simply an expected log-likelihood ratio, and the Poisson model has a multiscale likelihood factorization involving binomial conditional probabilities, we find (adopting a dyadic analysis and ignoring the coarsest scale term) that

$$\begin{aligned}
K(p_{\boldsymbol{\theta}}, p_{\boldsymbol{\theta}'}) &\sim \sum_{j,k} E\left[\log \frac{p(X_{j+1,2k}|X_{j,k}, \omega_{j,k})}{p(X_{j+1,2k}|X_{j,k}, \omega'_{j,k})}\right] \\
&= \sum_{j,k} \theta_{j,k}\left[\omega_{j,k} \log\left(\frac{\omega_{j,k}}{\omega'_{j,k}}\right) + (1-\omega_{j,k}) \log\left(\frac{1-\omega_{j,k}}{1-\omega'_{j,k}}\right)\right]. 
\end{aligned} \tag{34}$$

Next, define a hypercube (actually we use just the boundary or shell) in $\boldsymbol{\omega}$-space in analogy to that in (30) by specifying (i) $\theta_{0,0}$ is fixed and known, (ii) $\omega_{j,k} \equiv 1/2 \ \forall j \neq j^*$, and (iii) $\omega_{j,k} = 1/2 + s_k \Delta$ for $j = j^*$, where $s_k = \pm 1$ is an unknown sign and $\Delta \equiv \Delta(j^*) \in (0, 1/2)$ is a perturbation of known magnitude. Then, for this particular subproblem, estimation of $\boldsymbol{\theta}$ reduces to estimation of the $s_k$. Since, for appropriately defined values of $j^*$ and $\Delta$, our hypercube will be contained within the set of $\boldsymbol{\theta}$ induced by our function space $\Theta = \mathrm{BV}(C)$, it follows using (34) that

$$\frac{1}{N} \sum_{k=0}^{2^{j^*}-1} \theta_{j^*,k} \inf_{\hat{s}_k} \sup_{s_k} r(\hat{s}_k, s_k) \tag{35}$$

lower-bounds the quantity in (33), where

$$r(\hat{s}, s) \equiv E\left[\left(\frac{1}{2} + s\Delta\right) \log\left(\frac{1/2 + s\Delta}{1/2 + \hat{s}\Delta}\right) + \left(\frac{1}{2} - s\Delta\right) \log\left(\frac{1/2 - s\Delta}{1/2 - \hat{s}\Delta}\right)\right]. \tag{36}$$

However, the optimization problem $\inf \sup r(\hat{s}, s)$ is equivalent to a standard decision problem with binomial observations, a two-point action space, and 0/1 loss. Therefore, restricting attention to estimators of the form $\hat{s} = \pm 1$ we find that the expression for $r(\hat{s}, s)$ can be simplified to

$$(37) \qquad \left[ 2\Delta \log \left( \frac{1/2 + \Delta}{1/2 - \Delta} \right) \right] \Pr(\hat{s} \neq s).$$

Neyman–Pearson theory then yields that each individual probability $\Pr(\hat{s}_k \neq s_k)$ is minimized by the estimator $\hat{s}_k = \text{sgn}(X_{j^*+1,2k} - X_{j^*+1,2k+1})$, and therefore approximately equal to

$$(38) \qquad p^* \equiv \Pr\left\{ Z < -(\theta_{j^*+1,2k} - \theta_{j^*+1,2k+1})/\sqrt{\theta_{j^*,k}} \right\}$$

for sufficiently large $\theta_{j^*+1,2k}$ and $\theta_{j^*+1,2k+1}$, where $Z$ is a standard normal random variable. Note that by construction of our hypercube the value $p^*$ is the same for all $k$.

Finally, we address the issue of choice of $\Delta$ and $j^*$. First note that functions within our hypercube are simply piecewise constant functions of a fixed magnitude. Recalling the definition of the $\theta_{j,k}$ by average-sampling and the relation $\theta_{j+1,2k} = \theta_{j,k}\omega_{j,k}$, simple calculations yield that the height of any drop or rise in this function is simply proportional to $2\Delta$ [where the constant of proportionality is due to $\int_0^1 \theta(t)\,dt$, which may be arbitrarily set to 1]. With respect to the total variation seminorm defined in (13), this implies that for such functions to be in the space $\text{BV}(C)$ we must have $\Delta \leq (C/4)2^{-j^*}$. Next, analogous to the Gaussian case, we choose $j^*$ to satisfy the constraint $2^{3j/2} = (C/2)N^{1/2}$. This choice can be motivated, for example, by selecting that value of $j$ for which the signal-to-noise ratio in the empirical Haar coefficients $\langle \mathbf{X}, h_{j,k} \rangle$ is 1.

Combining the expressions in (35), (37) and (38), and exploiting the fact that the $\theta_{j^*,k}$ are equal for all $k$, we find that our lower bound on the minimax risk in (33) behaves like

$$(39) \qquad \frac{2^{j^*}}{N} \theta_{j^*,0} \left[ 2\Delta \log \left( \frac{1/2 + \Delta}{1/2 - \Delta} \right) \right] p^*.$$

To complete our proof, we note that the argument of the probability in (38) is in fact the Haar coefficient signal-to-noise ratio (i.e., the ratio of the expected value to the standard deviation), and hence $p^* \approx \Pr(Z < -1)$ may be treated as a constant in (39). Additionally, since $\theta_{j,k} = \theta_{0,0}/2^j$ and $\theta_{0,0} = N \int_0^1 \theta(t)\,dt$, it follows that $\frac{2^{j^*}}{N}\theta_{j^*,0}$ reduces to a constant. Lastly, a Taylor series expansion shows that the term within brackets behaves like $8\Delta^2$. From our definition of $\Delta$ and choice of $j^*$ it follows that $\Delta \sim N^{-1/3}$, from which the minimax lower bound rate of $N^{-2/3}$ stated in the theorem follows.

For the multinomial case, the proof proceeds in an almost identical manner. As the multinomial likelihood too has a multiscale likelihood factorization with binomial conditional probabilities, the expression in (34) holds in this case as well, although with $\theta_{j,k}$ replaced by $n\theta_{j,k}$. The hypercube is defined as before, but with $\theta_{0,0} \equiv 1$ by virtue of $\theta(\cdot)$ being a density, and the same underlying decision problem and optimal solution result. The discussion leading to our choice of $\Delta$ remains unchanged, and the constraint of $2^{3j/2} \sim N^{1/2}$ (with appropriately defined constants) again renders $p^*$ a constant. Finally, the quantity to the left of the bracketed expression in (39) now is $(2^{j^*}/N) n\theta_{j^*,0}$, which is approximately 1 by the defining condition $n \sim N$ in our specification of model (M). The rest of the argument is identical to that of the Poisson case. □

PROOF OF THEOREM 6—MODELS (P) AND (M). As the proof of Theorem 4 for models (P) and (M) exploited the $L_2$ approximation error properties of optimal nonlinear Haar approximants in a manner analogous to that of model (G), proof of Theorem 6 and the other results for these models follow using precisely the same modifications described earlier for model (G). □

**6. Discussion.** In this paper we lay a succinct conceptual foundation for the existence of certain multiscale likelihood factorizations. We also establish adaptivity and near-optimality of certain multiscale complexity penalized likelihood estimators, based on these factorizations, through study of their risk behavior. A key feature of our formulation and analysis is that canonical models of continuous, count and categorical data—Gaussian, Poisson and multinomial—are handled with common estimators, algorithms and risk analysis. The properties of our estimators derive essentially from their ability to exploit the fact that the decoupling inherent in the underlying multiscale factorizations for these models mirrors the decomposition deriving from an associated Haar basis.

In a sense this paper can be viewed also as providing some degree of explanation of and justification for the performance of other earlier work by the authors and colleagues with multiscale factorizations in specific methodological contexts, such as the analysis of Poisson time series [Kolaczyk (1999a, b)] and images [Timmermann and Nowak (1999)], Poisson linear inverse problems [Nowak and Kolaczyk (2000)] and the spatial analysis of continuous and count data in geography [Kolaczyk and Huang (2001)]. The multiscale likelihood approach analyzed here, based on average-sampling of the continuous object, fits quite naturally in many imaging applications in which the instrumentation involves spatially binning photon detections and a Poisson model. Similar comments apply in histogram and density estimation contexts involving binned data and the multinomial model. Moreover, it is in contexts like that of this last paper that some particularly interesting aspects of the flexibility of the multiscale likelihood framework come to light.

For example, in the geographical analysis of census data the notion of "multiscale" might arise through a desire to consider the effects of various levels of geopolitical aggregation (e.g., towns, counties, states) on, say, changes in population dynamics between two decennial censuses. A spatial–temporal analysis of this type may be set up and executed using a framework completely analogous to that underlying the recursive dyadic partitioning estimator considered in this paper. See Kolaczyk and Huang (2001) for details.

Lastly, we include two comments regarding details in our method of proof for risk analysis. First, we note the role that our adaptation of results of Li (1999) and Li and Barron (2000) plays in our derivation of upper bounds for the risk. The bound presented in Theorem 7 is quite general, and our subsequent usage of that bound suggests its usefulness in other problems of complexity penalized likelihood estimation in the nonparametric context. Second, we point out that, while our proof of the upper bounds on the risk (i.e., Theorem 4) does not explicitly use the decoupled structure of our likelihood factorizations, this structure does play a key role in our adaptation of the method of hyperrectangles for providing lower bounds on the risk (i.e., Theorem 5) in the Poisson and multinomial cases, wherein Kullback–Leibler divergences are simply $\ell_2$ risks in the standard Gaussian sequence model.

## APPENDIX

PROOF OF THEOREM 3. By virtue of the likelihood factorization (4), and the additivity of penalty function (8) in the multiscale indexing $I \in \mathcal{I}_{\mathrm{NT}}(\mathcal{P}^*)$, for all three estimators (9), (10) and (11) the objective function to be optimized in (7) is simply of the form $\sum_{I \in \mathcal{I}_{\mathrm{NT}}(\mathcal{P}^*)} W_I$, where the $W_I$ are functions of the data $\mathbf{X}$ and parameters $\boldsymbol{\theta}$. In particular, for each estimator there is at each $I$ a choice being made between $H_I^{(0)} : \{\omega_I \text{ is trivial}\}$, versus $H_I^{(1)} : \{\omega_I \text{ is nontrivial}\}$. The $W_I$ are either the negative log-likelihood under the null $H_I^{(0)}$ or that under the alternative $H_I^{(1)}$ plus a penalty in the amount of $2\lambda$.

In the case of thresholding, determining the value for $W_I$ for each $I$ corresponds to performing $N$ independent generalized likelihood ratio tests, and hence is of $O(N)$ algorithmic complexity trivially.

The case of recursive dyadic partition estimators follows arguments parallel to those in Donoho (1997). Specifically, any RDP of $[0,1)$ can be matched in one-to-one correspondence with a dyadic Haar basis in which there is a "hereditary" constraint on which coefficients are "kept" and which are "killed." That is, if a coefficient is to be kept, all of its ancestors must be kept as well; conversely, if a coefficient is killed, all of its descendents are killed as well. Hence, searching for an optimal RDP, say $\hat{\mathcal{P}} \preceq \mathcal{P}^*_{\mathrm{Dy}}$, is equivalent to a type of constrained thresholding. The constraint may be enforced by recursively moving from fine scales to coarse, and at each interval $I_{j,k}$ choosing the optimal subpartition $\hat{\mathcal{P}}(I_{j,k})$ on $I_{j,k}$ to be either (i) the union of the optimal subpartitions $\hat{\mathcal{P}}(I_{j+1,2k})$ and $\hat{\mathcal{P}}(I_{j+1,2k+1})$ on the children of $I_{j,k}$ or (ii) the trivial subpartition in which $I_{j,k}$ is partitioned no further. These decisions are based on the same type of generalized likelihood ratio tests underlying the thresholding estimators, indexed in the $I_{j,k} \in \mathcal{I}_{\mathrm{NT}}(\mathcal{P}^*_{\mathrm{Dy}})$. Therefore, this is the same as the optimal pruning algorithm for CART, which requires on the order of $N$ operations [Donoho (1997) and Breiman, Friedman, Olshen and Stone (1984), algorithm 10.1, page 294].

Finally, we consider the case of the general recursive partitioning estimators. First note there are a total of $\frac{N(N+1)}{2}$ unique subintervals $I \subseteq [0,1)$ that may be composed of the $N$ finest resolution intervals $I_i$. For any length $N_I \geq 1$, there are exactly $N - N_I + 1$ of these subintervals of length $N_I$. Also note that any interval of cardinality $N_I = m$ may be partitioned into two children intervals in exactly $m-1$ ways. Therefore, in total, among the $(N-1)!$ possible C-RPs, there are only $\sum_{m=1}^{N} m(N-m) = \frac{N^2(N+1)}{2} - \frac{N(N+1)(2N+1)}{6} \sim \frac{N^3}{6}$ unique parent–child pairs. By exploiting both this redundancy and the inheritance property underlying the RDP case, an efficient dynamic programming algorithm can be obtained. Here, however, the end result of the algorithm is not only an optimal partition $\hat{\mathcal{P}}$, but an accompanying C-RP $\mathcal{P}^*(\hat{\mathcal{P}})$ as well. Beginning with intervals $I$ of cardinality $N_I = 2$, and working recursively over $N_I = 3, 4, \ldots$, one can compute $W_I$ under both hypotheses $H_I^{(0)}$ and $H_I^{(1)}$, and pass the optimal decision (i.e., partition or not) for each interval "upward" (i.e., toward the coarsest interval $[0,1)$). That is, for a given interval $I$ of length $N_I = m$ [which may or may not appear in the definition of the final C-RP $\mathcal{P}^*(\hat{\mathcal{P}})$], we determine and record the optimal subpartition $\hat{\mathcal{P}}(I)$ on $I$ and the associated optimal sub-C-RP $\mathcal{P}^*(\hat{\mathcal{P}}(I))$, and we record the complexity value

$$\sum_{I' \in \mathcal{I}_{\mathrm{NT}}(\mathcal{P}^*(\hat{\mathcal{P}}(I)))} W_{I'}. \tag{40}$$

This optimal (sub)partition and its complexity value can be found via a maximization over $O(m)$ terms involving the corresponding optimal subpartitions and complexity values determined previously for each of the $m-1$ pairs of possible children of $I$. The maximization over all possible blocks in all possible C-RPs therefore requires roughly $N^3/6$ comparisons. Once we reach the top, we are left with $\hat{\mathcal{P}}$ and $\mathcal{P}^*(\hat{\mathcal{P}})$. □

**Acknowledgments.** The authors thank Andrew Barron and David Donoho for helpful discussions. They also thank the referees, Associate Editors and Editors involved at various stages for their comments and suggestions leading up to the final version of this paper.

# REFERENCES

Bar-Lev, S. K. and Enis, P. (1986). Reproducibility and natural exponential families with power variance functions. *Ann. Statist.* **14** 1507–1522. MR868315

Barndorff-Nielsen, O. (1978). *Information and Exponential Families in Statistical Theory.* Wiley, New York. MR489333

Barron, A., Birgé, L. and Massart, P. (1999). Risk bounds for model selection via penalization. *Probab. Theory Related Fields* **113** 301–413. MR1679028

Barron, A. R. and Cover, T. M. (1991). Minimum complexity density estimation. *IEEE Trans. Inform. Theory* **37** 1034–1054.

Breiman, L., Friedman, J., Olshen, R. and Stone, C. J. (1984). *Classification and Regression Trees.* Wadsworth, Belmont, CA. MR1111806

Daubechies, I. (1992). *Ten Lectures on Wavelets.* SIAM, Philadelphia. MR1162107

DeVore, R. A. (1998). Nonlinear approximation. In *Acta Numerica* **7** 51–150. Cambridge Univ. Press. MR1689432

Donoho, D. L. (1993). Unconditional bases are optimal bases for data compression and for statistical estimation. *Appl. Comput. Harmon. Anal.* **1** 100–115. MR1256530

Donoho, D. L. (1997). CART and best-ortho-basis: A connection. *Ann. Statist.* **25** 1870–1911. MR1474073

Donoho, D. L., Johnstone, I. M., Kerkyacharian, G. and Picard, D. (1995). Wavelet shrinkage: Asymptopia? (with discussion). *J. Roy. Statist. Soc. Ser. B* **57** 301–369. MR1323344

Donoho, D. L., Liu, R. and MacGibbon, B. (1990). Minimax risk over hyperrectangles, and implications. *Ann. Statist.* **18** 1416–1437. MR1062717

Girardi, M. and Sweldens, W. (1997). A new class of unbalanced Haar wavelets that form an unconditional basis for $L_p$ on general measure spaces. *J. Fourier Anal. Appl.* **3** 457–474. MR1468375

Joshi, S. W. and Patil, G. P. (1970). A class of statistical models for multiple counts. In *Random Counts in Scientific Work* (G. P. Patil, ed.) **2** 189–203. Pennsylvania State Univ. Press. MR287599

Kolaczyk, E. D. (1999a). Bayesian multiscale models for Poisson processes. *J. Amer. Statist. Assoc.* **94** 920–933. MR1723303

Kolaczyk, E. D. (1999b). Some observations on the tractability of certain multi-scale models. In *Bayesian Inference in Wavelet-Based Models. Lecture Notes in Statist.* **141** 51–66. Springer, New York. MR1699876

Kolaczyk, E. D. and Huang, H. (2001). Multiscale statistical models for hierarchical spatial aggregation. *Geographical Analysis* **33** 95–118.

Lauritzen, S. L. (1996). *Graphical Models.* Oxford Univ. Press. MR1419991

Li, Q. J. (1999). Estimation of mixture models. Ph.D. dissertation, Dept. Statistics, Yale Univ.

Li, Q. J. and Barron, A. R. (2000). Mixture density estimation. In *Advances in Neural Information Processing Systems* **12** 279–285. MIT Press, Cambridge, MA.

Nowak, R. D. (1999). Multiscale hidden Markov models for Bayesian image analysis. In *Bayesian Inference in Wavelet-Based Models. Lecture Notes in Statist.* **141** 243–265. Springer, New York. MR1699845

Nowak, R. D. and Kolaczyk, E. D. (2000). A statistical multiscale framework for Poisson inverse problems. *IEEE Trans. Inform. Theory* **46** 1811–1825. MR1790322

Sweldens, W. (1998). The lifting scheme: A construction of second generation wavelets. *SIAM J. Math. Anal.* **29** 511–546. MR1616507

Timmermann, K. E. and Nowak, R. D. (1999). Multiscale modeling and estimation of Poisson processes with application to photon-limited imaging. *IEEE Trans. Inform. Theory* **45** 846–862. MR1682515

Wilks, S. S. (1962). *Mathematical Statistics.* Wiley, New York. MR144404

Department of Mathematics
 and Statistics
Boston University
Boston, Massachusetts 02215
USA
e-mail: kolaczyk@math.bu.edu

Department of Electrical
 and Computer Engineering
University of Wisconsin
Madison, Wisconsin 53706
USA
e-mail: nowak@engr.wisc.edu